\documentclass[11pt]{article}

\usepackage{bussproofs}

\usepackage{pxfonts}
\usepackage{microtype}

\usepackage{natbib}

\begin{document}

\title{Comment on Mark Textor: `Brentano's Positing Theory of Existence'}
\author{Nils K\"urbis}
\date{}
\maketitle

\noindent \textsc{Note added April 2023.} This is the text of a commentary on a talk delivered by Mark Textor at King's College London in December 2015 (published version \citep{textorbrentano}). The text remains unchanged, except for corrections of typos and the layout. I have since changed my mind on a number of issues presented here, but it contains a few thoughts and an idea for a system of natural deduction that incorporates Textor's account of Brentano's positing theory of existence, as presented at the conference, which may not be devoid of interest. I have recently picked up on the idea in as yet unpublished work.

\section{The Trouble with Existence}
The standard introduction and elimination rules for quantifiers in free logic, with $\exists ! t$ standing for `$t$ exists', are:

\begin{center}
\bottomAlignProof
\AxiomC{$[\exists ! a]^i$}
\noLine
\UnaryInfC{$\Pi$}
\noLine
\UnaryInfC{$A_a^x$}
\RightLabel{$_i$}
\LeftLabel{$(\forall I)$ \ }
\UnaryInfC{$\forall x A$}
\DisplayProof\qquad
\bottomAlignProof
\AxiomC{$\forall x A$}
\AxiomC{$\exists  ! t$}
\LeftLabel{$(\forall E)$ \ }
\BinaryInfC{$A_t^x$}
\DisplayProof\bigskip

\bottomAlignProof
\AxiomC{$A_t^x$}
\AxiomC{$\exists ! t$}
\LeftLabel{$(\exists I) \ $}
\BinaryInfC{$\exists x A$}
\DisplayProof\qquad
\bottomAlignProof
\AxiomC{$\exists x A$}
\AxiomC{$[A_a^x]^i \ [\exists ! a]^j$}
\noLine
\UnaryInfC{$\Pi$}
\noLine 
\UnaryInfC{$C$}
\RightLabel{$_{i, j}$}
\LeftLabel{$(\exists E)$ \ }
\BinaryInfC{$C$}
\DisplayProof
\end{center}

\noindent where in $(\forall I)$, $a$ is not free in any assumptions on which $A$ depends except $\exists ! a$, and in $(\exists E)$ $a$ is not free in the minor premise $C$ and any formula it depends on, except $\exists ! a$ and $A_a^x$. 

These rules are unsatisfactory. The intended interpretation of $\exists ! t$ is `$t$ exists', but the rules don't tell us that. They are just as good if you read $\exists !$ as your favourite one-place predicate other than `exists'. The domain of quantification is then restricted to those things falling under that predicate, whatever it may be. So far, the rules are just an exercise in restricted quantification. 

The rules use $\exists !$, so they determine the meanings of the quantifiers only if $\exists !$ already has a meaning. We need a prior explanation of the meaning of $\exists !$ for the rules to impart the intended meaning onto the quantifiers. The question is how to do this in terms of rules of inference. 

We could try to add identity to solve the problem. Its elimination rule is the indiscernibility of identities. For the introduction rule we have two options. If $t=t$ is true whether $t$ refers or not, we use self-identity as an axiom. If $t=t$ is true only if $t$ refers, we use the restricted axiom $\forall x \ x=x$ or require $\exists ! t$ as a premise before we can deduce $t=t$: 

\begin{center}
$(=I_1)$ 
\AxiomC{}
\UnaryInfC{$t=t$}
\DisplayProof\quad
$(=I_2)$
\AxiomC{}
\UnaryInfC{$\forall x \ x=x$}
\DisplayProof\quad
$(=I_3)$
\AxiomC{$\exists ! t$}
\UnaryInfC{$t = t$}
\DisplayProof\quad
$(= E)$ \AxiomC{$t=u$}
\AxiomC{$A_t^x$}
\BinaryInfC{$A_u^x$}
\DisplayProof
\end{center}

\noindent Whichever introduction rule we use, $\exists ! t$ is now interderivable with $\exists x \ x=t$: 

\begin{center}
\AxiomC{}
\UnaryInfC{$t=t$}
\AxiomC{$\exists ! t$}
\BinaryInfC{$\exists x \ x=t$}
\DisplayProof\qquad
\AxiomC{}
\UnaryInfC{$\forall x \ x=x$}
\AxiomC{$\exists ! t$}
\BinaryInfC{$t=t$}
\AxiomC{$\exists !t$}
\BinaryInfC{$\exists x \ x=t$}
\DisplayProof\bigskip

\AxiomC{$\exists ! t$}
\UnaryInfC{$t=t$}
\AxiomC{$\exists !t$}
\BinaryInfC{$\exists x \ x=t$}
\DisplayProof\qquad
\AxiomC{$\exists x \ x=t$}
\AxiomC{$[a=t]^1$}
\AxiomC{$[\exists ! a]^1$}
\BinaryInfC{$\exists ! t$}
\RightLabel{$_1$}
\BinaryInfC{$\exists ! t$}
\DisplayProof
\end{center}

\noindent The universal quantifier occurs in $(=I_2)$, so it can only be understood if $\forall $ is understood. But $\forall$ in turn can only be understood if $\exists !$ is understood. That's circular and so won't help us understanding $\exists !$ or give a reason to exclude deviant interpretations of it. We should go for $(= I_1)$ or $(=I_3)$. As $(=I_2)$ and $(=I_3)$ are interderivable, it makes no difference to derivability which one we choose.

We could try the following line of thought. Identity is not restricted to a domain of things we are interested in for idiosyncratic reasons. It applies to everything. In the presence of identity, we can exclude deviant interpretations of $\exists !$, as $\exists ! t$ and $\exists x \ x=t$ are interderivable. For this line of thought to get off the ground, we cannot adopt $(=I_3)$, as it makes the meaning of $=$ dependent on the meaning of $\exists !$ (and therefore cannot exclude deviant interpretations either). 

We are left with $(=I_1)$. Sadly, it doesn't fare much better. Our line of thought makes the meaning of $\exists !$ dependent on the meaning of $=$, without making this explicit by the rules, but that cannot hide circularity. Even setting aside worries that we should be able to give the meaning of the existence predicate independently of identity, the dependence is problematic. We might as well \emph{define} $\exists ! t$ in terms of $=$, in the way suggested by its interderivability with $\exists x \ x=t$. To guarantee that $\exists !$ gets its intended meaning, we cannot understand $\exists !$ without understanding $=$. However, tying $\exists !$ to $=$ in turn appeals to $\exists$ (either in the proof of the equivalence of $\exists ! t$ with $\exists x \ x=t$ or in the definition of $\exists !$ in terms of $=$), but the meaning of $\exists $ is specified by rules appealing to $\exists !$. 

Neil Tennant proposes to add a rule of \emph{atomic denotation}, which goes some way towards alleviating the difficulties: An atomic proposition can only be true if all its singular terms refer. To secure the freedom of our logic, as $t=t$ is atomic, we have to reject $(= I_1)$. Tennant weakens it and requires an atomic proposition as premise from to infer $t=t$. Where $Ft$ is an atomic formula: 

\begin{center}
$(AD)$
\AxiomC{$Ft$}
\UnaryInfC{$\exists ! t$}
\DisplayProof\qquad
$(= I_4)$
\AxiomC{$Ft$}
\UnaryInfC{$t=t$}
\DisplayProof
\end{center} 

\noindent As $\exists ! t$ is atomic, it is interderivable with $t=t$ and also $\exists x \ x=t$. Contrary to the above, here the interderivability doesn't strike me as problematic.\footnote{Milne sees this as a problem in (\cite{milneexistence}), which is a discussion of Tennant's views on existence, reference and free logic, esp. (\cite{tennantabstraction}). Tennant's response (\cite{tennantmilne}) contains the rules given here.} $(AD)$ is a genuine introduction rule and there is nothing circular about it. The meaning of $\exists !$ and $=$ are now both dependent on a prior understanding of the meanings of some atomic propositions. But that's the case for all rules of inference. A more serious problem is that it is questionable whether the rule of atomic denotation in fact imparts on $\exists ! t$ the meaning of `$t$ exists', rather than `"$t$" denotes'. And the existence predicate, as Mark points out, should be distinguished from metalinguistic concepts like denotation. 

A formal aside: As $(\forall E)$ and $(\exists I)$ are effectively elimination rules for $\exists !$, deductions can contain maximal formulas of the form $\exists ! t$ that cannot be removed, namely when $\exists ! t$ is derived from an atomic formula and then used in an application of $(\forall E)$ or $(\exists I)$. Maybe that's fine: as the quantifier rules allow formulas of the form $\exists ! t$ to disappear, we only have a restricted subformula property, where some occurrence of $\exists !$ are discounted. Maybe we don't want to say that the meaning of $\exists !$ is determined \emph{exclusively} by the rules governing it, so that normalisability is an excessive requirement.

\section{Enter Brentano and Textor} 
The account free logic gives us of existence is not very satisfactory. What we need is some independent explanation of the concept of existence, ideally one that can be formally implemented in a logic. Mark, following Brentano, gives us a promising way of filling a logical gap. Mark suggests to revise Brentano's expressivist view of `existence': existence is a property that everything has. That's clearly what free logic must accept, too: in the rules for the quantifiers and $(AD)$, $\exists !$ is of exactly the same category and plays the same role as first-order predicates. Without a corresponding property, free logic wouldn't make sense! Nonetheless, Mark explains how we can make good use of Brentano's ideas: we can explain existence in terms of a primitive notion of acknowledgement and its dual `rejection', which are kinds of speech acts. `$t$ exists' means that if one acknowledges $t$, one acknowledges truly and if one rejects $t$, one rejects falsely. Mark adopts Rumfitt's idea that the meanings of the logical constants are given in terms of primitive speech acts of assertion $+$ and denial $-$ \citep{rumfittyesno} and proposes to extend it to the definition of the meaning of the existence predicate in terms of rules of inference that appeal to acknowledging and rejecting. For negation we have: 

\begin{center}
$(+ \neg I)$ \AxiomC{$- \ A$}
\UnaryInfC{$+ \ \neg A$}
\DisplayProof\quad
$(+\neg E)$ \AxiomC{$+\ \neg A$}
\UnaryInfC{$- \ A$}
\DisplayProof\qquad
$(-\neg I)$ \AxiomC{$- \ A$}
\UnaryInfC{$+\ \neg  A$}
\DisplayProof\quad
$(-\neg E)$ \AxiomC{$-\ \neg A$}
\UnaryInfC{$+ \ A$}
\DisplayProof
\end{center}

\noindent Mark proposes analogous rules for $\exists !$, where I adjust the notation for coherence with the rules above. $!$ stands for the acknowledging force and $/$ for the rejecting force: 

\begin{center}
$(\exists ! I_1)$ \AxiomC{$! \ t$}
\UnaryInfC{$+ \ \exists ! t$}
\DisplayProof\quad
$(\exists ! E_1)$ \AxiomC{$+\ \exists ! t$}
\UnaryInfC{$! \ t$}
\DisplayProof\qquad
$(\exists ! I_2)$ \AxiomC{$/  \ t$}
\UnaryInfC{$+ \ \neg \exists ! t$}
\DisplayProof\quad
$(\exists ! E_2)$ \AxiomC{$+\ \neg \exists ! t$}
\UnaryInfC{$/ \ t$}
\DisplayProof
\end{center}

\noindent Two minor points. First, as $\neg$ is used in the rules, the meaning of $\exists !$ is dependent on the meaning of $\neg$. So I couldn't understand $\exists !$ unless I understand $\neg$. But we can easily reformulate the rules by using the primitive force of denial in $(\exists ! I_2)$ and $(\exists ! E_2)$ instead: 

\begin{center}
$(\exists ! I_2')$ \AxiomC{$/  \ t$}
\UnaryInfC{$- \ \exists ! t$}
\DisplayProof\quad
$(\exists ! E_2')$ \AxiomC{$- \ \exists ! t$}
\UnaryInfC{$/ \ t$}
\DisplayProof
\end{center}

\noindent Second, Mark needs rules that specify that $/$ and $!$ are incompatible speech acts and what to do when we've arrived at such an impasse. The following suggest themselves: 

\begin{center}
\AxiomC{$!\ t$}
\AxiomC{$/ \ t$}
\BinaryInfC{$\bot$}
\DisplayProof\qquad
\AxiomC{$[+ \ \exists ! t]^i$}
\noLine
\UnaryInfC{$\Pi$}
\noLine
\UnaryInfC{$\bot$}
\RightLabel{$_i$}
\UnaryInfC{$/ \ t$}
\DisplayProof\qquad
\AxiomC{$[- \ \exists ! t]^i$}
\noLine
\UnaryInfC{$\Pi$}
\noLine
\UnaryInfC{$\bot$}
\RightLabel{$_i$}
\UnaryInfC{$! \ t$}
\DisplayProof
\end{center}

\noindent We can reformulate the rules for free logic in the light of Mark's discussion. The rules above are good only for the assertive use of the quantifiers. We also need rules for denial. Instead of appealing to $\exists !$, we can appeal directly to $! \ t$ in the rules in which $\exists ! t$ appears as a premise. Mark points out that, as speech acts cannot be embedded in sentences, but `$t$ exists' can, we need an existence predicate, just as we need $\neg$ in addition to denial. Maybe we can argue that we cannot eliminate $\exists !$ from the rules in which it occurs as a discharged hypothesis: it makes sense to discharge the assertion or denial of propositions, but it does not make sense to discharge acknowledgement or rejection of things. This would give another reason why we need an existence predicate. We arrive at the following rules, where $\alpha$ ranges over assertions and denials: 

\begin{center}
\bottomAlignProof
\AxiomC{$[+ \ \exists ! a]^i$}
\noLine
\UnaryInfC{$\Pi$}
\noLine
\UnaryInfC{$+ \ A$}
\RightLabel{$_i$}
\LeftLabel{$(+ \forall I) \ $}
\UnaryInfC{$+ \ \forall x A_x^a$}
\DisplayProof\qquad
\bottomAlignProof
\AxiomC{$+ \ \forall x A$}
\AxiomC{$! \ t$}
\LeftLabel{$(+ \forall E) \ $}
\BinaryInfC{$+ \ A_t^x$}
\DisplayProof\bigskip

\bottomAlignProof
\AxiomC{$+ \ A_t^x$}
\AxiomC{$ ! \ t$}
\LeftLabel{$(+ \exists I)$}
\BinaryInfC{$+ \ \exists x A$}
\DisplayProof\qquad
\bottomAlignProof
\AxiomC{$+ \ \exists x A$}
\AxiomC{$[+ \ A_a^x]^i \ [+ \ \exists ! a]^j$}
\noLine
\UnaryInfC{$\Pi$}
\noLine 
\UnaryInfC{$\alpha$}
\RightLabel{$_{i, j}$}
\LeftLabel{$(+ \exists E)$}
\BinaryInfC{$\alpha$}
\DisplayProof
\end{center}

\begin{center}
\bottomAlignProof
\AxiomC{$- \ A_t^x$}
\AxiomC{$ ! \ t$}
\LeftLabel{$(- \forall I)$}
\BinaryInfC{$- \ \forall x A$}
\DisplayProof\qquad
\bottomAlignProof
\AxiomC{$- \ \forall x A$}
\AxiomC{$[- \ A_a^x]^i \ [+ \ \exists ! a]^j$}
\noLine
\UnaryInfC{$\Pi$}
\noLine 
\UnaryInfC{$\alpha$}
\RightLabel{$_{i, j}$}
\LeftLabel{$(- \forall E)$}
\BinaryInfC{$\alpha$}
\DisplayProof\bigskip

\bottomAlignProof
\AxiomC{$[+ \ \exists ! a]^i$}
\noLine
\UnaryInfC{$\Pi$}
\noLine
\UnaryInfC{$- \ A$}
\RightLabel{$_i$}
\LeftLabel{$(- \exists I)$}
\UnaryInfC{$- \ \exists x A_x^a$}
\DisplayProof\qquad
\bottomAlignProof
\AxiomC{$- \ \exists x A$}
\AxiomC{$! \ t$}
\LeftLabel{$(- \exists E)$}
\BinaryInfC{$- \ A_t^x$}
\DisplayProof
\end{center}

\noindent These rules neither make use of $/ \ t$ nor of $-\ \exists ! t$. We could add a version of $(+ \ \forall I)$ that allows us to derive $/ \ t$ from $+\ \forall x A$ and $- \ A_t^x$, and a version of $(- \ \exists E)$ that allows us to derive $/ \ t$ from $- \ \exists xA$ and $+ \ A_t^x$, but that seems rather forced. I haven't suggested rules that allow the discharge of $- \ \exists ! t$, and I wonder what they could add. On the whole, at this point $/$ seems redundant. 

There remains some work to be done of how to integrate Brentano's and Mark's account of the meaning of $\exists !$ in a formal system, but it promises to be one that does not suffer from the circularities of the rules of free logic, and it is an account of the existence predicate, rather than denotation. Concerning identity, I suggest $(=I_2)$ is still best avoided, as we can probably make sense of identity independently of universal quantification, but Mark's account can go with either $(=I_1)$ or a version of $(=I_3)$, with $! \ t$ instead of $\exists ! t$ as the premise.

\section{Questions and Suggestions}
1. Assertion and denial, rejection and acknowledgement are all speech acts. What is the relation between assertion and acknowledgement on the one hand, and denial and rejection on the other? Are they species of a wider genus, so that there is \emph{one} speech act behind the use of $+$ and $!$ and another behind $-$ and $/$, only that it is applied to propositions in the first case, and to objects in the second? One reason why I have formulated $(+ \forall E)$, $(+\exists I)$, $(-\forall I)$ and $(-\exists E)$ by appealing directly to $!$ instead of $\exists !$ is that doing so demonstrates how they can all be integrated into a logic. Sticking to the old rules for the quantifiers makes Mark's rules for $!$ look rather like add-ons that are somewhat disconnected from the logic. If there are only two speech acts that can be applied equally to propositions and objects, that would give a neat unified account of the meanings of the logical constants and existence and possibly more, as suggested in the next point. Notice also that as $- \ A$ is equivalent to $+\ \neg A$, bilateralists give a fairly elaborate account of why we need both speech acts. There is nothing corresponding to $/$ that is equivalent to a combination of $!$ and $\neg$, as we cannot apply negation to objects. So there is an immediate need for $/$ that is absent for $-$.\bigskip

\noindent 2. I can acknowledge the existence of the the inventor of the zip fastener, but reject the existence of the present king of France. We can apply $!$ and $/$ also to complex terms, such as those constructed with the definite description operator $\iota$. We should have, e.g, $! \ \iota x Fx\vdash + \ F(\iota xFx)$, and our logic should then have two interlocking parts. One specifies the meanings of sentential connectives and quantifiers and of complex propositions on the basis of the conditions for the assertion and denial of atomic ones. Another part specifies the meanings of term forming operators and the conditions of acknowledgement and rejection of complex terms in terms of simpler ones. (We cannot reconstruct $/$ as $! \ (\iota x \ \neg \ x=t)$, which is correct only if there is only one object and it is not $t$!)\bigskip 

\noindent 3. Should there be something corresponding to the rule of atomic denotation? Can we, for instance, infer from the assertion of an atomic proposition $Ft$ that $t$ is acknowledged, and conversely, that if $t$ is rejected, an atomic proposition $Ft$ must be, too? If so, we have $+ \ Ft\vdash \ ! \ t$ and $/ \ t \vdash - \ Ft$.

\bigskip
\setlength{\bibsep}{0pt}
\bibliographystyle{chicago}
\bibliography{KurbisCommentOnTextor}
\end{document}